\documentclass[preprint,5p]{elsarticle}
\makeatletter
\def\ps@pprintTitle{%
 \let\@oddhead\@empty
 \let\@evenhead\@empty
 \def\@oddfoot{}%
 \let\@evenfoot\@oddfoot}
\makeatother

\usepackage{amsmath,amssymb,amsfonts,amsthm}
\usepackage{algorithmic}
\usepackage[boxed]{algorithm2e}
\usepackage{mathtools}
\usepackage{pgfplots}
\usepackage{tikz}
\usepackage{subfig}
\usepackage{euscript}
\usepackage{fontawesome}
\usepackage{pstricks}
\usepackage{comment}
\usepackage{enumitem}
\usepackage{esdiff}
\usepackage{mathabx} % FOR ARROWS
\usepackage{booktabs} % For thick hlines
\usetikzlibrary{decorations.pathreplacing}

\usepackage{pifont}
\pgfplotsset{compat=1.13}
\usepgfplotslibrary{fillbetween}
\usetikzlibrary{patterns}

\pgfplotsset{%
  colormap={bluewhite}{color(0cm)=(white);
  color(1cm)=(blue!50!black)}
}

\tikzset{
  font={\fontsize{8pt}{12}\selectfont}}

\pgfmathdeclarefunction{gauss}{2}{%
  \pgfmathparse{1/(#2*sqrt(2*pi))*exp(-((x-#1)^2)/(2*#2^2))}%
}


\newtheorem{theorem}{Theorem}
\newtheorem{proposition}{Proposition}

\newtheorem{definition}{Definition}

\newtheorem{remark}{Remark}

\SetKwInput{KwData}{Data}
\SetKwInput{KwReturn}{return}

\newcommand{\shrink}{\vspace{-0.3cm}}

\newcommand{\maximize}[1]{\underset{{#1}}{\text{max}}}
\newcommand{\st}{\text{s.t.}}
\newcommand{\set}{\EuScript}
\usepackage{mathtools}
\DeclarePairedDelimiterX{\norm}[1]{\lVert}{\rVert}{#1}
\newcommand*\circled[1]{\tikz[baseline=(char.base)]{
            \node[shape=circle,draw,inner sep=0.5pt] (char) {\footnotesize{#1}};}}
\newcommand\mydots{\hbox to 0.75em{.\hss.\hss.}}
\definecolor{red}{rgb}{0.4, 0.6, 0.8}
\definecolor{blue}{rgb}{0.16, 0.32, 0.75}

\usepackage{lineno,hyperref}
\modulolinenumbers[5]

\begin{document}
\begin{frontmatter}
\title{Electricity Market Equilibrium under Information Asymmetry}
\author{Vladimir Dvorkin Jr., Jalal Kazempour, Pierre Pinson}
\address{Technical University of Denmark, Elektrovej 325, 2800 Kongens Lyngby, Denmark}
\begin{abstract}
We study a competitive electricity market equilibrium with two trading stages, day-ahead and real-time. The welfare of each market agent is exposed to uncertainty (here from renewable energy production), while agent information on the probability distribution of this uncertainty is not identical at the day-ahead stage. We show a high sensitivity of the equilibrium solution to the level of information asymmetry and demonstrate economic, operational, and computational value for the system stemming from potential information sharing.
\end{abstract}
\end{frontmatter}
%\linenumbers

\section{Introduction}
With increasing shares of renewable energy resources, electricity markets are exposed to uncertainty associated with intermittent power supply. To accommodate this uncertainty, electricity trading in short term has been  arranged in several subsequent trading floors, such as day-ahead and real-time markets. At the first stage (day-ahead), market agents, such as power producers and consumers, compete to contract energy considering the forecast of renewable power production, while at the second stage (real-time) they settle power imbalances caused by forecast errors. 

Such competition can be modeled as a stochastic {\it equilibrium} problem in the sense of \cite{arrow1954existence}, where each market agent is a self-optimizer and maximizes its welfare, e.g. expected profit for producers or expected utility for consumers. Each agent computes a day-ahead decision while anticipating the real-time outcomes using its {\it private} information (e.g., probabilistic forecast) about uncertainty. To find the equilibrium among such agents, distributed algorithms as in \cite{gerard2018risk} have been proposed to let agents integrate their private forecasts. After a finite number of iterations, the agents find a set of equilibrium prices that reflect their private information and support the equilibrium. This distributed solution is promising for the operation of local markets \cite{8301552} or parts of larger systems \cite{Molzahn_2017}.

In contrast, large electricity markets, such as NordPool or CAISO, use a {\it centralized} optimization to compute the equilibrium, where the central entity called market operator collects bids from agents and clears the market on its own. To efficiently operate markets with renewables, it has been proposed to cast the centralized optimization as a two-stage stochastic model \cite{morales2012pricing,pritchard2010single}. This model considers that the operator generates a set of plausible renewable outcomes based on its {\it own} information about the probability distribution of renewable energy production, and clears the day-ahead market while accounting for the anticipated real-time imbalances.

The important property of the two problems is that they are equivalent in the case where all market agents in the equilibrium model optimize against the same probability distribution as that of the market operator in the centralized model. Under this scenario, the two problems yield the same market-clearing results and the centralized market is complete as it satisfies the preferences of the agents in the equilibrium problem. However, the equivalence between the centralized and equilibrium problems no longer holds when agents in the equilibrium problem optimize against different distributions. In this situation, the centralized market settlement is inefficient as it does not support the true preferences of agents.  We refer to this situation as {\it information asymmetry} that typically holds for many reasons. Naturally, market agents use different data and forecast tools to build uncertainty distributions. Furthermore, for a given set of plausible outcomes, agents may explicitly assign different probabilities depending on whether they are rational or not, for example, in the sense of prospect theory \cite{tversky1992advances}.

In this line, this letter analyzes a competition among electricity market agents that have asymmetric information about a common source of uncertainty. We propose an equilibrium model in which agents may assign different probabilities over a common set of renewable power production outcomes. We show that the centralized model in \cite{morales2012pricing, pritchard2010single} satisfies the preferences of all market agents only if they all agree on the probability distribution of renewable power production. We discuss the existence and uniqueness of the solution to the equilibrium problem, and refer the stability theory to discuss the challenges associated with its computation. With our analytic results, we point out a high sensitivity of equilibrium prices to the level of information asymmetry. We then propose a distributed algorithm to numerically assess the equilibrium outcomes. We eventually demonstrate that the system overall benefits from information sharing, as we show that for any asymmetry in agents information, there exist a loss of social welfare, an increase in real-time imbalances, and a decrease of the convergence rate of the distributed algorithm.

The letter is outlined as follows. In Section 2 we describe the setup and introduce the stochastic equilibrium model. We further proceed with an analytic solution to equilibrium prices as a function of the private information of agents in Section 3. In Section 4 we describe the distributed algorithm to compute equilibrium and provide extensive numerical experiments. All proofs are gathered in an Appendix.

\section{Problem statement}
\subsection{Main notation and assumptions}
We consider a finite set of uncertainty outcomes $\Omega$ indexed by $\omega=\{1,\dots,\Omega\}$. $\xi_{\omega}$ is the renewable power output that corresponds to outcome $\omega$. The renewable producers are not modeled as market agents, but represented as an aggregated stochastic in-feed. The controllable generation (consumption) is represented by a single producer (consumer). The dispatch of power producer at the day-ahead stage is denoted by $p\in\set{O}$, and it can be adjusted by $r_{\omega}\in\set{O}$ in real-time if outcome $\omega$ realizes. The set $\set{O}$ denotes the feasible operating region of the producer based on its technical constraints. The cost function of the producer is quadratic given by $c(x) = \frac{1}{2}\alpha x^2$, where $\alpha$ is a positive constant. The consumer procures energy at the day-ahead stage in amount of $d\in\set{K}$ that is subject to adjustment $l_{\omega}\in\set{K}$ in real-time if outcome $\omega$ realizes. The set $\set{K}$ exhibits the feasible region of the decision-making problem of the consumer. The utility of the consumer is described by concave function $u(x) = \gamma x - \frac{1}{2}\beta x^2$, where $\gamma$ and $\beta$ are positive constants. We assume that both $\set{O}$ and $\set{K}$ are convex and compact sets. The \textcolor{black}{dual} price in scenario $\omega$ is denoted by $\lambda_{\omega}$ in the optimization problem. \textcolor{black}{Its counterpart in the equilibrium problem} is denoted by $\tilde{\lambda}_{\omega}$. \textcolor{black}{In this work, we do not consider network, subsidies and unit commitment constraints, that often cause negative electricity prices \cite{deng2015cost}, and exclusively focus on perfect competition. Therefore, both $\lambda_{\omega}$ and $\tilde{\lambda}_{\omega}$} belong to a compact set of non-negative reals $\Lambda_{+}$.

\subsection{Centralized model for market-clearing problem}
Consider a centralized market organization, where the market operator collects bids of agents and finds  socially optimal contracts $\{p,d\}$ at the day-ahead  stage, followed by real-time recourse decisions $\{r_{\omega},l_{\omega}\}_{\forall \omega}$. The market operator integrates its own information about underlying uncertainty  that is described by a finite set of probabilities $\{\pi_{\omega}^{\text{mo}}\}_{\forall\omega}$ assigned to uncertain outcomes. This yields
\begin{subequations}
\begingroup
\allowdisplaybreaks
\label{stoch_problem}
\begin{align}
    \maximize {p,r_{\omega},d,l_{\omega}} & \big[u(d)-c(p)\big] + \sum_{\omega \in \Omega}\pi_{\omega}^{\text{mo}}\big[u(l_{\omega})-c(r_{\omega})\big], \label{total_cost}\\
    \st \quad & p + r_{\omega} + \xi_{\omega} - d - l_{\omega} \geq 0 : \lambda_{\omega}, \quad \forall \omega \in \Omega,\label{PB_stoch}\\
    & (p,r_{\omega}) \in \set{O}, \; (d,l_{\omega}) \in \set{K}, \quad \forall \omega \in \Omega,
\end{align}
\endgroup
\end{subequations}
where objective function \eqref{total_cost} represents the expected social welfare seen by the market operator, and constraint \eqref{PB_stoch} enforces the power balance for each outcome of renewable energy production. A set of dual prices $\{\lambda_{\omega}\}_{\forall \omega}$ shows the sensitivity of the expected social welfare to the stochastic in-feed and, therefore, is an implicit function of the information of market operator. Hence, the outcomes for market participants are subject to the information available to the market operator.
\textcolor{black}{
\begin{remark}
As dual prices fall into non-negative domain, \eqref{PB_stoch} is cast as an inequality constraint that is binding in optimum. If one allows for negative prices, \eqref{PB_stoch} has to be specified as an equality constraint. 
\end{remark}
\begin{remark}
The real-time electricity price in outcome $\omega$ anticipated by the market operator at the day-ahead stage is the probability-removed price $\frac{\lambda_{\omega}}{\pi_{\omega}^{\text{mo}}}$ \cite{morales2012pricing}. 
\end{remark}
}
\begin{remark}
Unlike settings in \cite{morales2012pricing,pritchard2010single}, we do not explicitly model the day-ahead power balance constraint. Instead, we use the notion of price convergence between day-ahead and real-time stages \cite{hogan2016virtual} to obtain the day-ahead \textcolor{black}{electricity} price as \textcolor{black}{$\lambda^{\text{DA}} = \Sigma_{\omega}\pi_{\omega}^{\text{mo}}\frac{\lambda_{\omega}}{\pi_{\omega}^{\text{mo}}}=\Sigma_{\omega}\lambda_{\omega}$}.
\end{remark}

% We highlight that unlike market models in \cite{morales2012pricing,pritchard2010single}, we do not explicitly model the day-ahead power balance constraint, and use the notion of price convergence between day-ahead and real-time stages \cite{hogan2016virtual} to obtain the day-ahead price as $\lambda^{\text{DA}} = \Sigma_{\omega}\lambda_{\omega}$.

\subsection{Equilibrium model for market-clearing problem}
We now introduce an equilibrium model given by a set of individual optimization of three agents, i.e.,
\begin{subequations}\label{equilibrium_model}
\begin{align}
    \maximize{\tilde{\lambda}_{\omega}\in\Lambda_{+}}\;J_{\omega}^{\text{ps}}&:=-\tilde{\lambda}_{\omega}\big[p + r_{\omega} + \xi_{\omega} - d - l_{\omega}\big], \; \forall \omega\in\Omega,\label{price_setter} \\
    \maximize{(p,r_{\omega}) \in \set{O}} J^{\text{p}}&:=\sum_{\omega \in \Omega}\pi_{\omega}^{\text{p}}\left[\frac{\tilde{\lambda}_{\omega}}{\pi_{\omega}^{\text{p}}}(p + r_{\omega})-c(r_{\omega})\right]-c(p), \label{profit_max}\\
    \maximize{(d,l_{\omega}) \in \set{K}} J^{\text{c}}&:= \sum_{\omega\in\Omega}\pi_{\omega}^{\text{c}} \left[u(l_{\omega}) - \frac{\tilde{\lambda}_{\omega}}{\pi_{\omega}^{\text{c}}}(d + l_{\omega})\right] + u(d), \label{utility_max}
\end{align}
\end{subequations}

The price-setting agent solves \eqref{price_setter} and optimizes a set of \textcolor{black}{equilibrium} prices $\{\tilde{\lambda}_{\omega}\}_{\forall\omega}$ in response to the value of the system imbalance for each outcome of renewable production. For any surplus of generation, problem \eqref{price_setter} yields zero price, while it yields a strictly positive price in case of generation shortage. The power producer optimizes its first- and second-stage decisions $p$ and $\{r_{\omega}\}_{\forall \omega}$ in \eqref{profit_max} to maximize the expected profit for a given set of prices $\{\tilde{\lambda}_{\omega}\}_{\forall\omega}$. In its optimization, the producer integrates its own information about the uncertain in-feed characterized by a finite set of probabilities $\{\pi_{\omega}^{\text{p}}\}_{\forall\omega}$. Finally, the consumer computes optimal first-stage and recourse decisions $d$ and $\{l_{\omega}\}_{\forall\omega}$ in \eqref{utility_max} to maximize its expected utility using its own information set $\{\pi_{\omega}^{\text{c}}\}_{\forall\omega}$. Observe, that agents in \eqref{profit_max} and \eqref{utility_max} use the probability-removed prices obtained by dividing the equilibrium prices by the associated probabilities \cite{morales2012pricing}. The probability-removed prices define the actual electricity price that each agent expects to receive once uncertainty is resolved.

The three problems are interconnected in the sense that the problem of the price-setter is parametrized by the decisions of the producer and consumer, while their problems are conditioned by the price provided by the price-setting agent. \textcolor{black}{Similarly to the centralized problem (1a), equilibrium prices provide the sensitivity of the expected social welfare with respect to the marginal change in random in-feed.} Therefore, a set of equilibrium prices $\{\tilde{\lambda}_{\omega}\}_{\forall\omega}$ is implicitly a function of the information that agents integrate into their optimization problems. 
\begin{proposition}\label{Prop_ex_uniq}
The solution to the equilibrium problem \eqref{equilibrium_model} exists and is unique for any agent information sets.
\end{proposition}
% \begin{proof}
% See \ref{Proof_ex_uniq}.
% \end{proof}
\begin{remark}
The proof of Proposition \ref{Prop_ex_uniq} relies on the strict monotonicity of agent preferences. In the case of linear preferences, other approaches would be required (see \cite[Chapter 2]{facchinei2007finite}).
\end{remark}

\subsection{Relation between centralized and equilibrium models}
The equivalence between the centralized and equilibrium market-clearing models is established with the following proposition.
\begin{proposition}\label{Prop_equiv}
Let $\pi_{\omega}^{\text{mo}} = \pi_{\omega}^{p} = \pi_{\omega}^{c}, \forall \omega \in \Omega.$ Then, there exists a set of prices $\{\tilde{\lambda}_{\omega}^{\star}\}_{\forall\omega}$ that yields the optimal solution  $p^{\star},d^{\star},\{r_{\omega}^{\star},l_{\omega}^{\star}\}_{\forall\omega}$ in the equilibrium model \eqref{equilibrium_model} that solves the centralized model \eqref{stoch_problem}. Moreover, $\tilde{\lambda}_{\omega}^{\star}=\lambda_{\omega}^{\star}, \forall \omega.$
\end{proposition}
% \begin{proof}
% See \ref{Proof_Prop_1}.
% \end{proof}
However, this equivalence no longer holds when the information of market agents about the renewable in-feed in the equilibrium model is different from that of the market operator in the centralized model. In this scenario, the prices in \eqref{stoch_problem} and \eqref{equilibrium_model} are not necessarily identical as they depend on different information sets, making the market based on \eqref{stoch_problem} incomplete in terms of information. In the following we study model \eqref{equilibrium_model} that reveals the true equilibrium state among agents with private information on uncertainty. Eventually, we show that the system overall benefits when agents agree on a common information set that completes the market.

\section{Analytic solution for equilibrium prices}
Let us define the demand excess function for renewable power outcome $\omega$ as
$z_{\omega} = d + l_{\omega} - p - r_{\omega} - \xi_{\omega}.$
We derive the optimality conditions associated with \eqref{profit_max} and \eqref{utility_max} to define variables $d,l_{\omega}, p,$  and $r_{\omega}$ as a function of equilibrium prices $\boldsymbol{\tilde{\lambda}}$. Assuming the agent constraints are not binding, the demand excess function writes as:
%
%
%
% that can be expressed as the following function of prices considering the optimality conditions of \eqref{profit_max} and \eqref{utility_max} and assuming that constraints on agents' decisions are not binding:
\begin{align*}
    z_{\omega}(\boldsymbol{\tilde{\lambda}}) = \frac{\gamma-\Sigma_{\omega}\tilde{\lambda}_{\omega}}{\beta} + \frac{\pi_{\omega}^{\text{c}}\gamma-\tilde{\lambda}_{\omega}}{\pi_{\omega}^{\text{c}}\beta} - \frac{\Sigma_{\omega}\tilde{\lambda}_{\omega}}{\alpha} - \frac{\tilde{\lambda}_{\omega}}{\pi_{\omega}^{\text{p}}\alpha} - \xi_{\omega}.
\end{align*}

\begin{figure}[]
\centering
\resizebox{9cm}{!}{%
\includegraphics[]{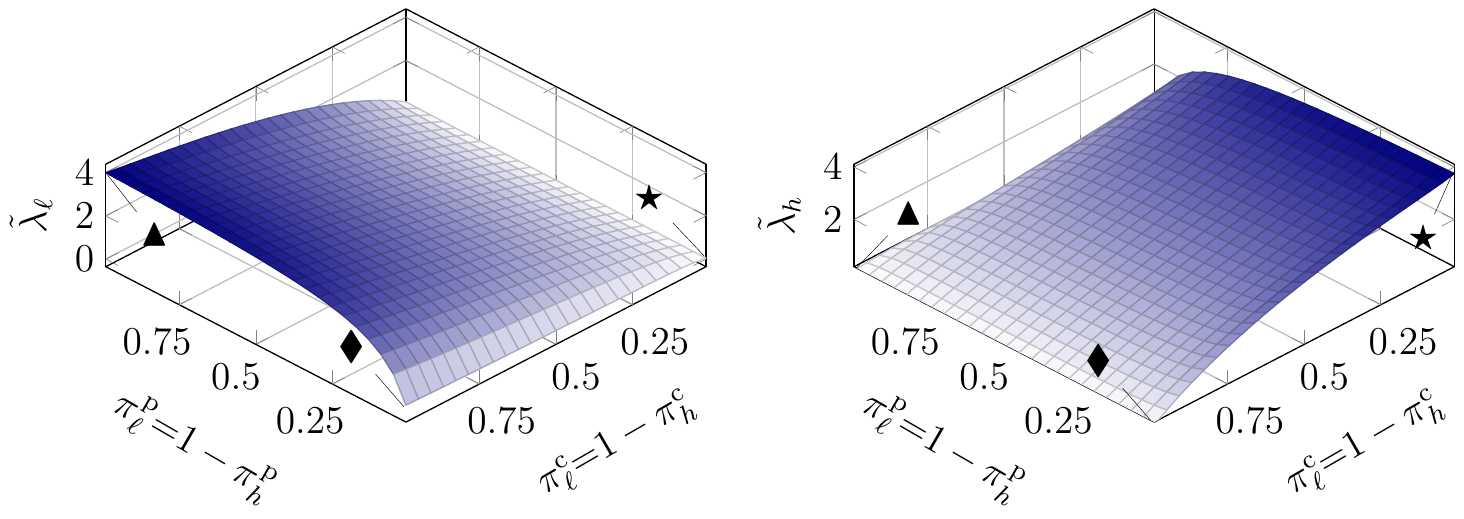}
}
\vspace{-0.6cm}
\caption{Equilibrium prices $\tilde{\lambda}_{\ell}$ and $\tilde{\lambda}_{h}$ as a function of probabilities that agents assign to the two uncertainty outcomes. The black markers indicate the three boundary equilibrium cases.}
\label{price_3d}
\vspace{-0.4cm}
\end{figure}

By solving $z_{\omega}(\boldsymbol{\tilde{\lambda}}) = 0, \forall\omega\in\Omega$, we obtain a closed-form characterization of equilibrium prices as a function of probabilities that agents assign to uncertain outcomes. In the interest of illustration, let us consider a set $\Omega\in\{\mathit{h},\ell\}$ with only two outcomes with $\xi_{\ell}=1$ and $\xi_{\mathit{h}}=3$. For any agent it holds that $\pi_{\ell}+\pi_{\mathit{h}}=1$. Let $\alpha=1.5$, $\beta=0.3$, and $\gamma=5$.
Figure \ref{price_3d} depicts the two equilibrium prices $\tilde{\lambda}_{\ell}$ and $\tilde{\lambda}_{h}$ as a function of $\pi_{\ell}$ and $\pi_{\mathit{h}}$.  We find a clear relationship between the equilibrium prices and agent information. For instance in case ($\blacktriangle$), when producer assigns the whole probability mass to outcome $\ell$, it leads to a nearly zero price associated with outcome $h$. A similar situation holds in the opposite case ($\bigstar$). In a quite critical case ($\blacklozenge$) with highly asymmetric assignment of probabilities, the equilibrium yields almost zero prices for both outcomes. Moreover, we find that the day-ahead price, i.e., $\tilde{\lambda}^{\text{DA}} = \tilde{\lambda}_{\ell} + \tilde{\lambda}_{h}$, attains its maximum value when both agents have symmetric information, i.e., $\pi_{\ell}^{\text{p}} = \pi_{\ell}^{\text{c}}$.

\begin{figure}[]
\centering
\subfloat[$\pi_{\ell}^{\text{p}}=0.5,\pi_{\ell}^{\text{c}}=0.5$]{\includegraphics[width= 1.5in]{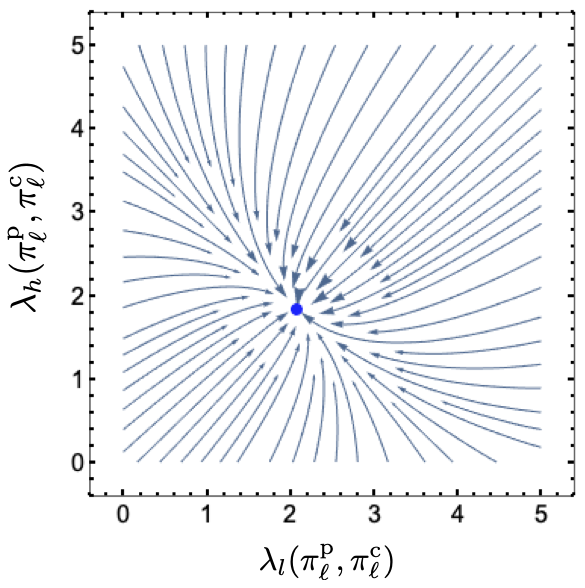}}\;
\subfloat[$\pi_{\ell}^{\text{p}}=0.99,\pi_{\ell}^{\text{c}}=0.5$]{\includegraphics[width= 1.5in]{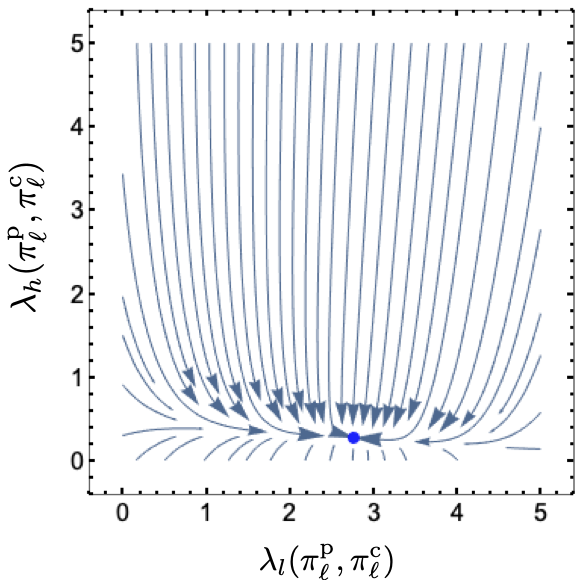}}\;
\vspace{-0.3cm}
\caption{Equilibrium point and vector field around equilibrium point in case of (a) symmetric and (b) asymmetric information.}
\label{three_eq}
\shrink
\end{figure}

As shown in \cite{gerard2018risk}, the unstable equilibrium may not be computable by standard distributed algorithms.  To verify the stability of the equilibrium solution under different assignments of probabilities, we consider a dynamic price adjustment process as the following first order differential equation \cite{arrow1958stability}:
\begin{align}\label{diff}
\frac{\text{d}\boldsymbol{\tilde{\lambda}}(t)}{\text{d}t} = \tau \boldsymbol{z}(\boldsymbol{\tilde{\lambda}}(t)),\quad \boldsymbol{\tilde{\lambda}}(0)=\boldsymbol{\tilde{\lambda}}_{0},
\end{align}
where $\tau$ is some positive constant, and $\boldsymbol{\tilde{\lambda}}_{0}$ is a vector of initial prices. We discuss the stability of the equilibrium solution using the following proposition.
\begin{proposition}[Adapted from \cite{bellman2008stability}]
If $\boldsymbol{\tilde{\lambda}}$ is a solution of \eqref{diff} and all the eigenvalues of the Jacobian matrix of $\boldsymbol{z}$ have strictly negative real parts, then $\boldsymbol{\tilde{\lambda}}$ is locally stable. If at least one eigenvalue has strictly positive real part, then $\boldsymbol{\tilde{\lambda}}$ is unstable.
\end{proposition}

By verifying the eigenvalues of the Jacobian of $\boldsymbol{z}$, we find that for any assignment of probabilities, the equilibrium solution is locally stable and, thus, supposedly computable. However, we observe that for asymmetric cases the ratio of the two eigenvalues significantly increases. This ratio heavily affects the convergence rate for gradient search algorithms \cite{goh2017why}, as illustrated by  vector fields in Figure \ref{three_eq} for some choice of $\boldsymbol{\tilde{\lambda}}_{0}$. In particular, in Figure \ref{three_eq}(a), a gradient search is almost uniform in both directions $\tilde{\lambda}_{\ell}$ and $\tilde{\lambda}_{h}$, while in Figure \ref{three_eq}(b) the gradient in direction $\tilde{\lambda}_{\ell}$ is notably smaller than that in direction $\tilde{\lambda}_{h}$. In the following section, we will demonstrate that the high eigenvalue ratios significantly affect the convergence of the distributed market-clearing algorithms.

\vspace{-0.4cm}
\section{Equilibrium computation}
In this section, we first introduce a distributed algorithm to compute the solution to the equilibrium problem. We then describe the setup and provide numerical results.
\subsection{Algorithm}
To compute the equilibrium solution, we use a distributed algorithm that naturally embodies the Walrasian tatonnement \cite{uzawa1960walras}. The price-setter problem \eqref{price_setter} updates the prices based on the optimal response of the producer and consumer optimization problems \eqref{profit_max} and \eqref{utility_max}, respectively.

We first show that the price-setter optimization \eqref{price_setter} reduces to a single analytic expression.
\begin{proposition}\label{Prop_price_setter_solution}
Consider the response of producer $p^{\nu},\{r_{\omega}^{\nu}\}_{\forall\omega}$ and the response of consumer $d^{\nu},\{l_{\omega}^{\nu}\}_{\forall\omega}$ to a set of prices $\{\tilde{\lambda}_{\omega}^{\nu-1}\}_{\forall\omega}$ at some iteration $\nu$.  Then, the solution of \eqref{price_setter} converges to optimum over iterations using
$$\tilde{\lambda}_{\omega}^{\nu} = \max \Big\{0,\tilde{\lambda}_{\omega}^{\nu-1} - \rho \big[p^{\nu} + r_{\omega}^{\nu} + \xi_{\omega} - d^{\nu} - l_{\omega}^{\nu}\big]\Big\}, \; \forall \omega \in \Omega,$$
for some positive constant $\rho$.
\end{proposition}

Using an analytic expression for the price-update, we can compute the solution of the equilibrium problem \eqref{equilibrium_model} using Algorithm 1. \textcolor{black}{As objective function of each agent is strictly monotone in decision variables and its feasibility set is convex and compact,  the algorithm provably converges to the global optimum for $\nu\rightarrow\infty$ with rate $\mathcal{O}(\frac{1}{\nu})$, given that the solution exists \cite{falsone2017dual}.} The algorithm is implemented \textcolor{black}{in JuMP environment \cite{dunning2017jump} in Julia, and the source code is available in the e-companion \cite{dvorkin2019}}.

\begin{algorithm}[h]
\DontPrintSemicolon
\SetAlgoLined
\KwData{$\nu_{\text{MAX}}$, $\rho$, $\tilde{\lambda}_{\omega}^{0} \; \forall \omega$, $\epsilon$}
 \For{$\nu$ from 1 to $\nu_{\text{MAX}}$}{
    \circled{\normalsize{1}} For $\{\tilde{\lambda}_{\omega}^{\nu-1}\}_{\forall\omega}$, update producer response
    $$p^{\nu},\{r_{\omega}^{\nu}\}_{\forall\omega}\leftarrow\underset{(p,r_{\omega}) \in \set{O}}{\text{argmax}}\quad   J^{\text{p}}(p,r_{\omega})$$
    \circled{\normalsize{2}} For $\{\tilde{\lambda}_{\omega}^{\nu-1}\}_{\forall\omega}$, update consumer response
    $$d^{\nu},\{l_{\omega}^{\nu}\}_{\forall\omega}\leftarrow\underset{(d,l_{\omega}) \in \set{K}}{\text{argmax}}\quad   J^{\text{c}}(d,l_{\omega})$$
    \circled{\normalsize{3}} For $p^{\nu},\{r_{\omega}^{\nu}\}_{\forall\omega}$ and $d^{\nu},\{l_{\omega}^{\nu}\}_{\forall\omega}$, update prices:
    $$\tilde{\lambda}_{\omega}^{\nu} = \max \Big\{0,\tilde{\lambda}_{\omega}^{\nu-1} - \rho \big[p^{\nu} + r_{\omega}^{\nu} + \xi_{\omega} - d^{\nu} - l_{\omega}^{\nu}\big]\Big\}$$
    \circled{\normalsize{4}} Return $\epsilon$-equilibrium prices and dispatch if:
    $$\norm[\big]{p^{\nu} + r_{\omega}^{\nu} + \xi_{\omega} - d^{\nu} - l_{\omega}^{\nu}}^{2} \leq \epsilon, \; \forall \omega\in\Omega,$$
    otherwise go to \circled{\normalsize{1}}.
 }
\caption{Solution algorithm}
\label{ALG}
\end{algorithm}
\vspace{-0.4cm}

\subsection{Setup}
We choose $\alpha = 1.5$, $\gamma=5$, $\beta=0.3$, $\{\tilde{\lambda}_{\omega}^{0}\}_{\forall\omega} = 0$, $\rho=\epsilon=10^{-5}$. The outcomes of uncertain renewable production are described by 100 samples drawn from a normal distribution $\mathcal{N}(\mu,\sigma^2)$ with $\mu=1.5$ and $\sigma^2=0.25$. \textcolor{black}{The rationale behind these parameters lies in the fact that the producer and consumer are willing to trade energy for any realization of wind power production, whereas the wind fluctuations bring about observable impacts on the market-clearing outcomes. The practical choice of these parameters is subject to the specifics of a given power system, e.g. cost/utility structure and wind penetration level}. 
        
We consider the reference distribution $\circled{\scriptsize{R}}$ that assigns equally likely probabilities over 100 samples. We then generate a series of distributions that tweak either mean or variance of the reference distribution $\circled{\scriptsize{R}}$ using the probability weighting function of the following form \cite[Eq.(3)]{gonzalez1999shape}:
\begin{align}
    \Phi(\boldsymbol{\xi}) = \frac{\delta\big[\Phi^{\circled{\scriptsize{\tiny{R}}}}(\boldsymbol{\xi})\big]^{\gamma}}{\delta\big[\Phi^{\circled{\scriptsize{\tiny{R}}}}(\boldsymbol{\xi})\big]^{\gamma} + \big[1-\Phi^{\circled{\scriptsize{\tiny{R}}}}(\boldsymbol{\xi})\big]^{\gamma}}, \label{PWF}
\end{align}
where $\Phi$ represents the cumulative distribution function of stochastic renewable production, $\delta\in\mathbb{R}_{+}$ primarily affects the mean of the reference distribution $\circled{\scriptsize{R}}$, and $\gamma\in\mathbb{R}_{+}$ primarily impacts the variance. By applying \eqref{PWF} to the reference distribution for different $\delta$ and $\gamma$, we obtain a collection of probability assignments to the same set of outcomes. Table \ref{dist_data_tab} summarizes the distributions that we use in the following analysis.

In our setup, consumer always optimizes against the reference distribution $\circled{\scriptsize{R}}$, while producer optimizes against one of the distributions in Table \ref{dist_data_tab}. When producer uses $\circled{\scriptsize{R}}$ in its local optimization, the equilibrium solution corresponds to the symmetric case, and any deviation from $\circled{\scriptsize{R}}$ corresponds to the asymmetric equilibrium.

\subsection{Numerical results}
\begin{table}
\footnotesize
\centering
\caption{Descriptive statistics of distributions: \protect\circled{\scriptsize{$\mu$}}-labeled distributions primarily tweak the mean of the reference distribution \protect\circled{\tiny{R}}, while \protect\circled{\scriptsize{$\sigma$}}-labeled distributions primarily tweak the variance of \protect\circled{\tiny{R}}.}
\begin{tabular}{c|ccc|c|ccc}
\specialrule{1pt}{1pt}{1pt}
Label & $\circled{$\mu$}_{3}^{\uparrow}$ & $\circled{$\mu$}_{2}^{\uparrow}$ & $\circled{$\mu$}_{1}^{\uparrow}$ & $\circled{\scriptsize{R}}$ & $\circled{$\mu$}_{1}^{\downarrow}$ & $\circled{$\mu$}_{2}^{\downarrow}$ & $\circled{$\mu$}_{3}^{\downarrow}$ \\
\hline
$\mu$ & 2.02 & 1.79 & 1.65 & 1.56 & 1.34 & 1.22 & 1.07  \\
$\sigma^2$ & 
\textcolor{black}{0.35} & 
\textcolor{black}{0.35} & 
\textcolor{black}{0.34} & 
\textcolor{black}{0.33} & 
\textcolor{black}{0.31} & 
\textcolor{black}{0.30} & 
\textcolor{black}{0.27} \\
\specialrule{1pt}{1pt}{1pt}
Label & $\circled{$\sigma$}_{3}^{\uparrow}$ & $\circled{$\sigma$}_{2}^{\uparrow}$ & $\circled{$\sigma$}_{1}^{\uparrow}$ & $\circled{\scriptsize{R}}$ & $\circled{$\sigma$}_{1}^{\downarrow}$ & $\circled{$\sigma$}_{2}^{\downarrow}$ & $\circled{$\sigma$}_{3}^{\downarrow}$ \\
\hline
$\mu$ & 1.63 & 1.60 & 1.57 & 1.56 & 1.55 & 1.55 & 1.56 \\
$\sigma^2$ & 
\textcolor{black}{1.62} & 
\textcolor{black}{0.92} & 
\textcolor{black}{0.54} & 
\textcolor{black}{0.33} & 
\textcolor{black}{0.10} & 
\textcolor{black}{0.04} & 
\textcolor{black}{0.02} \\
\specialrule{1pt}{1pt}{1pt}
\end{tabular}
\label{dist_data_tab}
\vspace{-0.6cm}
\end{table}

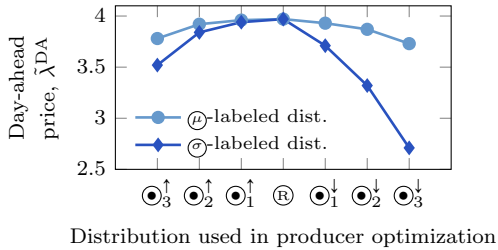
\begin{figure}
\centering
%\resizebox{0.5\textwidth}{!}{%
\begin{tikzpicture}
\begin{axis}[
xmin = 1,
xmax = 9,
ymin = 2.5,
ymax = 4.1,
xticklabels={$\circled{$\bullet$}_{3}^{\uparrow}$ , $\circled{$\bullet$}_{2}^{\uparrow}$ , $\circled{$\bullet$}_{1}^{\uparrow}$ , $\circled{\tiny{R}}$ , $\circled{$\bullet$}_{1}^{\downarrow}$ , $\circled{$\bullet$}_{2}^{\downarrow}$ , $\circled{$\bullet$}_{3}^{\downarrow}$},xtick={2,...,8},
typeset ticklabels with strut,  % <-- added
enlarge x limits=false,
height=3.75cm,
width=6cm,
ylabel style={align=center},
ylabel={Day-ahead \\ price, $\tilde{\lambda}^{\text{DA}}$},
legend pos=south west,
legend style={draw=none,fill=none},
xlabel={Distribution used in producer optimization}
]
\addplot[red, line width = 0.35mm, solid, mark=*,mark options={solid}] table [x=case, y=mean]{data/prices.dat};
\addlegendentry{\scriptsize{\circled{\tiny$\mu$}-labeled dist.}};
\addplot[blue, line width = 0.35mm, every mark/.append style={solid, fill=blue},mark=diamond*] table [x=case, y=var]{data/prices.dat};
\addlegendentry{\scriptsize{\circled{\tiny$\sigma$}-labeled dist.}};
\end{axis}
\end{tikzpicture}
%}
\shrink
\caption{Impacts of information asymmetry on the day-ahead price.}
\label{price}
\shrink
\end{figure}

\begin{figure}
\centering
%\resizebox{1\textwidth}{!}{%
\begin{tikzpicture}
\begin{axis}[
xshift=-4cm,
yshift=-4.5cm,
xmin = 1,
xmax = 9,
ymin = 2.25,
ymax = 4.5,
xticklabels={$\circled{\tiny$\mu$}_{3}^{\uparrow}$ , $\footnotesize\circled{\tiny$\mu$}_{2}^{\uparrow}$ , $\footnotesize\circled{\tiny$\mu$}_{1}^{\uparrow}$ , $\footnotesize\circled{\tiny{R}}$ , $\footnotesize\circled{\tiny$\mu$}_{1}^{\downarrow}$ , $\footnotesize\circled{\tiny$\mu$}_{2}^{\downarrow}$ , $\footnotesize\circled{\tiny$\mu$}_{3}^{\downarrow}$},
typeset ticklabels with strut,  % <-- added
enlarge x limits=false,
xtick={2,...,8},
height=3.85cm,
width=5cm,
ylabel style={align=center},
ylabel={Day-ahead \\ contracts},
]
% M+3	2.02	0.59	2.52	4.07
\addplot[red , line width = 0.35mm]  (1.75,4.07) -- (2.25,4.07);
\addplot[blue, line width = 0.35mm]  (1.75,2.52) -- (2.25,2.52);
\addplot[gray,line width = 0.35mm] (2,2.52) -- (2,4.07);
% M+2	1.79	0.59	2.61	3.59
\addplot[red , line width = 0.35mm]  (2.75,3.59) -- (3.25,3.59);
\addplot[blue, line width = 0.35mm]  (2.75,2.61) -- (3.25,2.61);
\addplot[gray,line width = 0.35mm] (3,2.61) -- (3,3.59);
% M+1	1.65	0.58	2.64	3.46
\addplot[red , line width = 0.35mm]  (3.75,3.46) -- (4.25,3.46);
\addplot[blue, line width = 0.35mm]  (3.75,2.64) -- (4.25,2.64);
\addplot[gray,line width = 0.35mm] (4,2.64) -- (4,3.46);
% M0	1.56	0.58	2.65	3.43
\addplot[red , line width = 0.35mm]  (4.75,3.43) -- (5.25,3.43);
\addplot[blue, line width = 0.35mm]  (4.75,2.65) -- (5.25,2.65);
\addplot[gray,line width = 0.35mm] (5,2.65) -- (5,3.43);
% M-1	1.34	0.56	2.62	3.56
\addplot[red , line width = 0.35mm]  (5.75,3.56) -- (6.25,3.56);
\addplot[blue, line width = 0.35mm]  (5.75,2.62) -- (6.25,2.62);
\addplot[gray,line width = 0.35mm] (6,2.62) -- (6,3.56);
% M-2	1.22	0.55	2.58	3.78
\addplot[red , line width = 0.35mm]  (6.75,3.78) -- (7.25,3.78);
\addplot[blue, line width = 0.35mm]  (6.75,2.58) -- (7.25,2.58);
\addplot[gray,line width = 0.35mm] (7,2.58) -- (7,3.78);
% M-3	1.07	0.52	2.49	4.22
\addplot[red , line width = 0.35mm]  (7.75,4.22) -- (8.25,4.22);
\addplot[blue, line width = 0.35mm]  (7.75,2.49) -- (8.25,2.49);
\addplot[gray,line width = 0.35mm] (8,2.49) -- (8,4.22);
\end{axis}
\begin{axis}[yshift=-4.5cm,
xmin = 1,
xmax = 9,
xticklabels={$\circled{\tiny$\sigma$}_{3}^{\uparrow}$ ,
$\circled{\tiny$\sigma$}_{2}^{\uparrow}$ ,
$\circled{\tiny$\sigma$}_{1}^{\uparrow}$ ,
$\circled{\tiny{R}}$ ,
$\circled{\tiny$\sigma$}_{1}^{\downarrow}$ ,
$\circled{\tiny$\sigma$}_{2}^{\downarrow}$ ,
$\circled{\tiny$\sigma$}_{3}^{\downarrow}$},
typeset ticklabels with strut,  % <-- added
enlarge x limits=false,
xtick={2,...,8},
height=3.85cm,
width=5cm,
xlabel={Distribution used in producer optimization},
xlabel style={xshift=-2.0cm}
]
% V+3	1.63	1.27	2.35	4.92
\addplot[red, line width = 0.35mm]  (1.75,4.92) -- (2.25,4.92);
\addplot[blue, line width = 0.35mm]  (1.75,2.35) -- (2.25,2.35);
\addplot[gray,line width = 0.35mm] (2,2.35) -- (2,4.92);
% V+2	1.6	    0.96	2.56	3.86
\addplot[red , line width = 0.35mm]  (2.75,3.86) -- (3.25,3.86);
\addplot[blue, line width = 0.35mm]  (2.75,2.56) -- (3.25,2.56);
\addplot[gray,line width = 0.35mm] (3,2.56) -- (3,3.86);
% V+1	1.57	0.73	2.63	3.52
\addplot[red , line width = 0.35mm]  (3.75,3.52) -- (4.25,3.52);
\addplot[blue, line width = 0.35mm]  (3.75,2.63) -- (4.25,2.63);
\addplot[gray,line width = 0.35mm] (4,2.63) -- (4,3.52);
% V0	1.56	0.58	2.65	3.43
\addplot[red , line width = 0.35mm]  (4.75,3.43) -- (5.25,3.43);
\addplot[blue, line width = 0.35mm]  (4.75,2.65) -- (5.25,2.65);
\addplot[gray,line width = 0.35mm] (5,2.65) -- (5,3.43);
% V-1	1.55	0.31	2.47	4.31
\addplot[red , line width = 0.35mm]  (5.75,4.31) -- (6.25,4.31);
\addplot[blue, line width = 0.35mm]  (5.75,2.47) -- (6.25,2.47);
\addplot[gray,line width = 0.35mm] (6,2.47) -- (6,4.31);
% V-2	1.55	0.21	2.21	5.61
\addplot[red , line width = 0.35mm]  (6.75,5.61) -- (7.25,5.61);
\addplot[blue, line width = 0.35mm]  (6.75,2.21) -- (7.25,2.21);
\addplot[gray,line width = 0.35mm] (7,2.21) -- (7,5.61);
% V-3	1.56	0.13	1.81	7.63
\addplot[red , line width = 0.35mm]  (7.75,7.63) -- (8.25,7.63);
\addplot[blue, line width = 0.35mm]  (7.75,1.81) -- (8.25,1.81);
\addplot[gray,line width = 0.35mm] (8,1.81) -- (8,7.63);
\end{axis}
%\draw[step=0.25,black,thin] (0.75,-3.25) grid (2.25,-2.25);
\draw[red, line width = 0.35mm] (0.8,-2.45) -- (1,-2.45) node[anchor=north, xshift = 0.15cm, yshift = 0.25cm] {\textcolor{black}{$d$}};
\draw[blue, line width = 0.35mm] (0.8,-3.1) -- (1,-3.1) node[anchor=north, xshift = 0.15cm, yshift = 0.25cm] {\textcolor{black}{$p$}};
\draw[gray, line width = 0.35mm] (0.9,-2.45) -- (0.9,-3.1);
\draw[gray] (0.75,-3.225) -- (0.75,-2.275) -- (2.35,-2.275) -- (2.35,-3.225) -- (0.75,-3.225);
\draw [decorate,decoration={brace,amplitude=4pt,mirror,raise=2pt},yshift=0pt,gray]
(0.8,-3.1) -- (0.8,-2.45) node [black,midway,xshift=0.9cm,yshift=0.05cm] {\footnotesize
\textcolor{black}{mismatch}};
\end{tikzpicture}
%}
\shrink
\caption{Impacts of information asymmetry on the contracted quantities of consumer ($d$) and producer ($p$) at the day-ahead stage.}
% \vspace{-0.5cm}
\label{dispatch}
\shrink
\end{figure}
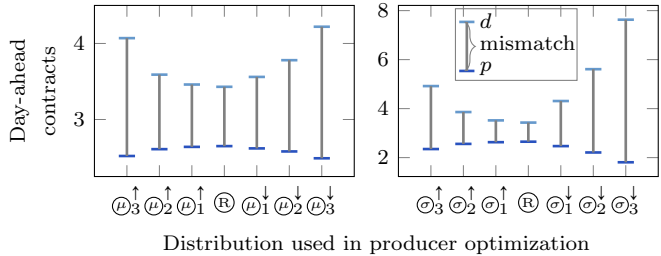

We first consider the impact of information asymmetry on the \textcolor{black}{electricity} price at the day-ahead stage depicted in Figure \ref{price}. We see that \textcolor{black}{it} is maximized when the two agents use the same uncertainty distribution $\circled{\scriptsize{R}}$. Any deviation from $\circled{\scriptsize{R}}$ in producer optimization decreases \textcolor{black}{the day-ahead} price. The resulting price-supported day-ahead contracts illustrated in Figure \ref{dispatch} show that such deviations in terms of either mean or variance lead to  increasing power mismatch between controllable generation and consumption.

Next, we compute the realization of the social welfare for each uncertainty outcome for the fixed day-ahead decisions of the producer and consumer. They are computed considering symmetric information and some asymmetric information cases. The results are summarized in Figure \ref{SW}.  We observe that the social welfare improves in larger realizations of renewable output, and records the maximum when producer employs $\circled{\scriptsize{R}}$. For any deviation of the producer from the reference distribution, we find a social loss, that is smaller for deviations in terms of the mean rather than variance for given distributions. Moreover, we see the welfare reduces more significantly if the producer assigns smaller variance relatively to that of the consumer.

Finally, we show how the computational performance of the algorithm is affected by the asymmetry of information. Table \ref{num_iter} collects the number of iterations required by the algorithm to converge along with the ratio between the largest and the smallest eigenvalues of Jacobian of the demand excess function. We see that apart from the case of $\circled{$\mu$}_{1}^{\uparrow}$ distribution, the asymmetry of agent information yields larger ratio of eigenvalues, and thus requires more iterations to converge. Moreover, for a highly asymmetric case of a low-variance distribution $\circled{$\sigma$}_{3}^{\downarrow}$, this ratio boosts so that the algorithm does not converge for any iteration limit. \textcolor{black}{The computational time of each iteration, though, is not affected by information asymmetry and kept below a few milliseconds for all distributions.}
\vspace{-0.5cm}

% \pgfplotsset{every tick label/.append style={font=\footnotesize}}

\definecolor{forestgreen}{rgb}{0.4, 0.6, 0.8}
\definecolor{wildwatermelon}{rgb}{0.16, 0.32, 0.75}

\begin{figure*}[t]
\resizebox{0.265\textwidth}{!}{%
\begin{tikzpicture}
    \begin{axis}[
    enlargelimits=0,
    xtick={1,20,40,60,80,100},
    xlabel={outcome \#},
    ylabel={Social welfare},
    xlabel near ticks,
    ylabel near ticks,
    legend pos=south east,
    legend style={draw=none,fill=none},
    legend columns=-1,
    height=3.75cm,
    width=5cm,
    ]
    \addplot[name path=d1,forestgreen!90,line width = 0.35mm]   table [x=x,y=base,col sep=space,mark=none] {data/SW_data.dat};
    \addlegendentry{$\circled{\tiny{R}}$}
    \addplot[name path=d01,wildwatermelon!90,line width = 0.35mm]   table [x=x,y=y1,col sep=space,mark=none] {data/SW_data.dat};
    \addlegendentry{$\circled{$\mu$}_{3}^{\uparrow}$}
    \addplot [
        thick,
        color=wildwatermelon,
        fill=wildwatermelon!30,
        fill opacity=0.5
    ]
    fill between[of=d1 and d01];
    \end{axis}
    \begin{axis}[
    xshift=4.5cm,
    enlargelimits=0,
    xtick={1,20,40,60,80,100},
    xlabel={outcome \#},
    xlabel near ticks,
    ylabel near ticks,
    legend pos=south east,
    legend style={draw=none,fill=none},
    legend columns=-1,
    height=3.75cm,
    width=5cm
    ]
    \addplot[name path=d1,forestgreen!90,line width = 0.35mm]   table [x=x,y=base,col sep=space,mark=none] {data/SW_data.dat};
    \addlegendentry{$\circled{\tiny{R}}$}
    \addplot[name path=d01,wildwatermelon!90,line width = 0.35mm]   table [x=x,y=y2,col sep=space,mark=none] {data/SW_data.dat};
    \addlegendentry{$\circled{$\mu$}_{3}^{\downarrow}$}
    \addplot [
        thick,
        color=wildwatermelon,
        fill=wildwatermelon!30,
        fill opacity=0.5
    ]
    fill between[of=d1 and d01];
    \end{axis}
    \begin{axis}[
    xshift=9cm,
    enlargelimits=0,
    xtick={1,20,40,60,80,100},
    xlabel={outcome \#},
    xlabel near ticks,
    ylabel near ticks,
    legend pos=south east,
    legend style={draw=none,fill=none},
    legend columns=-1,
    height=3.75cm,
    width=5cm
    ]
    \addplot[name path=d1,forestgreen!90,line width = 0.35mm]   table [x=x,y=base,col sep=space,mark=none] {data/SW_data.dat};
    \addlegendentry{$\circled{\tiny{R}}$}
    \addplot[name path=d01,wildwatermelon!90,line width = 0.35mm]   table [x=x,y=y3,col sep=space,mark=none] {data/SW_data.dat};
    \addlegendentry{$\circled{$\sigma$}_{3}^{\uparrow}$}
    \addplot [
        thick,
        color=wildwatermelon,
        fill=wildwatermelon!30,
        fill opacity=0.5
    ]
    fill between[of=d1 and d01];
    \end{axis}
    \begin{axis}[
    xshift=13.5cm,
    enlargelimits=0,
    xtick={1,20,40,60,80,100},
    xlabel={outcome \#},
    xlabel near ticks,
    ylabel near ticks,
    legend pos=south east,
    legend style={draw=none,fill=none},
    legend columns=-1,
    height=3.75cm,
    width=5cm
    ]
    \addplot[name path=d1,forestgreen!90,line width = 0.35mm]   table [x=x,y=base,col sep=space,mark=none] {data/SW_data.dat};
    \addlegendentry{$\circled{\tiny{R}}$}
    \addplot[name path=d01,wildwatermelon!90,line width = 0.35mm]   table [x=x,y=y4,col sep=space,mark=none] {data/SW_data.dat};
    \addlegendentry{$\circled{$\sigma$}_{3}^{\downarrow}$}
    \addplot [
        thick,
        color=wildwatermelon,
        fill=wildwatermelon!30,
        fill opacity=0.5
    ]
    fill between[of=d1 and d01];
    \end{axis}
\end{tikzpicture}
}
\vspace{-0.3cm}
\caption{Social welfare for each outcome of renewable production. \textcolor{black}{It is defined as a total system surplus for the fixed day-ahead decisions, i.e., $\text{SW}_{\omega}=\big[u(d^{\star})-c(p^{\star})\big] + \big[u(l_{\omega})-c(r_{\omega})\big]$.} The 100 outcomes are ordered from the smallest to largest. The colored area between the curves shows the welfare loss caused by asymmetry of information.}
\label{SW}
\shrink
\end{figure*}
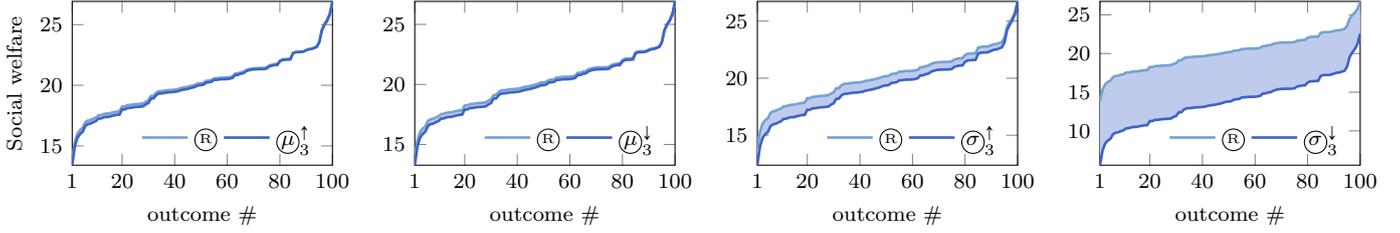

\begin{table}
\footnotesize
\centering
\label{tab_performance}
\caption{Number of iterations and the ratio between the largest and smallest eigenvalues for various distributions in producer optimization.}
\begin{tabular}{c|ccc|c|ccc}
\specialrule{1pt}{1pt}{1pt}
Label & $\circled{$\mu$}_{3}^{\uparrow}$ & $\circled{$\mu$}_{2}^{\uparrow}$ & $\circled{$\mu$}_{1}^{\uparrow}$ & $\circled{\scriptsize{R}}$ & $\circled{$\mu$}_{1}^{\downarrow}$ & $\circled{$\mu$}_{2}^{\downarrow}$ & $\circled{$\mu$}_{3}^{\downarrow}$ \\
\hline
\# iters & 369 & 347 & 312 & 321 & 362 & 372 & 381  \\
ratio & 2.4 & 2.2 & 2.0 & 2.0 & 2.2 & 2.3 & 2.5 \\
\specialrule{1pt}{1pt}{1pt}
Label & $\circled{$\sigma$}_{3}^{\uparrow}$ & $\circled{$\sigma$}_{2}^{\uparrow}$ & $\circled{$\sigma$}_{1}^{\uparrow}$ & $\circled{\scriptsize{R}}$ & $\circled{$\sigma$}_{1}^{\downarrow}$ & $\circled{$\sigma$}_{2}^{\downarrow}$ & $\circled{$\sigma$}_{3}^{\downarrow}$ \\
\hline
\# iters & 379 & 372 & 353 & 321 & 359 & 499 & $\infty$ \\
ratio & 2.7 & 2.4 & 2.2 & 2.0 & 18.8 & 1.8e3 & 1.8e7 \\
\specialrule{1pt}{1pt}{1pt}
\end{tabular}
\label{num_iter}
\shrink
\end{table}

\section*{References}
\bibliography{mybibfile}

\begin{thebibliography}{10}
\expandafter\ifx\csname url\endcsname\relax
  \def\url#1{\texttt{#1}}\fi
\expandafter\ifx\csname urlprefix\endcsname\relax\def\urlprefix{URL }\fi
\expandafter\ifx\csname href\endcsname\relax
  \def\href#1#2{#2} \def\path#1{#1}\fi

\bibitem{arrow1954existence}
K.~J. Arrow, G.~Debreu, Existence of an equilibrium for a competitive economy,
  Econometrica: Journal of the Econometric Society (1954) 265--290.

\bibitem{gerard2018risk}
H.~G{\'e}rard, V.~Lecl{\`e}re, A.~Philpott, On risk averse competitive
  equilibrium, Operations Research Letters 46~(1) (2018) 19--26.

\bibitem{8301552}
F.~Moret, P.~Pinson, Energy collectives: a community and fairness based
  approach to future electricity markets, IEEE Transactions on Power Systems,
  in print.

\bibitem{Molzahn_2017}
D.~K. Molzahn, F.~D{\"o}rfler, H.~Sandberg, S.~H. Low, S.~Chakrabarti,
  R.~Baldick, J.~Lavaei, A survey of distributed optimization and control
  algorithms for electric power systems, IEEE Transactions on Smart Grid 8~(6)
  (2017) 2941--2962.

\bibitem{morales2012pricing}
J.~M. Morales, A.~J. Conejo, K.~Liu, J.~Zhong, Pricing electricity in pools
  with wind producers, IEEE Transactions on Power Systems 27~(3) (2012)
  1366--1376.

\bibitem{pritchard2010single}
G.~Pritchard, G.~Zakeri, A.~Philpott, A single-settlement, energy-only electric
  power market for unpredictable and intermittent participants, Operations
  research 58~(4-part-2) (2010) 1210--1219.

\bibitem{tversky1992advances}
A.~Tversky, D.~Kahneman, Advances in prospect theory: Cumulative representation
  of uncertainty, Journal of Risk and uncertainty 5~(4) (1992) 297--323.

\bibitem{deng2015cost}
L.~Deng, B.~F. Hobbs, P.~Renson, What is the cost of negative bidding by wind?
  a unit commitment analysis of cost and emissions, IEEE Transactions on Power
  Systems 30~(4) (2015) 1805--1814.

\bibitem{hogan2016virtual}
W.~W. Hogan, Virtual bidding and electricity market design, The Electricity
  Journal 29~(5) (2016) 33--47.

\bibitem{facchinei2007finite}
F.~Facchinei, J.-S. Pang, Finite-dimensional variational inequalities and
  complementarity problems, Springer Science \& Business Media, 2007.

\bibitem{arrow1958stability}
K.~J. Arrow, L.~Hurwicz, On the stability of the competitive equilibrium,
  Econometrica: Journal of the Econometric Society (1958) 522--552.

\bibitem{bellman2008stability}
R.~Bellman, Stability theory of differential equations, Courier Corporation,
  2008.

\bibitem{goh2017why}
G.~Goh, \href{http://distill.pub/2017/momentum}{Why momentum really works}
  (2017).
\newline\urlprefix\url{http://distill.pub/2017/momentum}

\bibitem{uzawa1960walras}
H.~Uzawa, Walras' tatonnement in the theory of exchange, The Review of Economic
  Studies 27~(3) (1960) 182--194.

\bibitem{falsone2017dual}
A.~Falsone, K.~Margellos, S.~Garatti, M.~Prandini, Dual decomposition for
  multi-agent distributed optimization with coupling constraints, Automatica 84
  (2017) 149--158.

\bibitem{dunning2017jump}
I.~Dunning, J.~Huchette, M.~Lubin, {JuMP}: A modeling language for mathematical
  optimization, SIAM Review 59~(2) (2017) 295--320.

\bibitem{dvorkin2019}
V.~Dvorkin, J.~Kazempur, P.~Pinson,
  \href{https://github.com/wdvorkin/EMEuIA}{Electronic companion -
  {E}lectricity market equilibrium under information asymmetry}.
\newline\urlprefix\url{https://github.com/wdvorkin/EMEuIA}

\bibitem{gonzalez1999shape}
R.~Gonzalez, G.~Wu, On the shape of the probability weighting function,
  Cognitive psychology 38~(1) (1999) 129--166.

\end{thebibliography}
\vspace{-0.4cm}

\appendix
\section{Proof of Proposition \ref{Prop_ex_uniq}}\label{Proof_ex_uniq}
\subsection{Preliminaries}
We connect the solution of the equilibrium problem to the solution of variational inequalities.
\begin{definition}
Consider a mapping $F:\mathbb{R}^{n}\rightarrow\mathbb{R}^{n}$ and a set $K\subseteq\mathbb{R}$. A solution set $\mathrm{SOL}$($K$,$F$) to the variational inequality problem $\mathrm{VI}$($K$,$F$) is a vector $x^{\star} \in K$ such that $\langle F(x^{\star}),x-x^{\star}\rangle \geq 0, \; \forall x \in K$.
\end{definition}
We use the results from \cite{facchinei2007finite} to establish the existence and uniqueness of the equilibrium solution.
\begin{theorem}[Corollary 2.2.5 \cite{facchinei2007finite}] \label{theorem_ex}
Suppose that $K$ is a compact and convex set, and that the mapping $F$ is continuous. Then, the set $\mathrm{SOL}$($K$,$F$) is nonempty and compact.
\end{theorem}
\begin{theorem}[Theorem 1.3.1 \cite{facchinei2007finite}] \label{theorem_un}
Let $F: U\rightarrow\mathbb{R}$ be continuously differentiable on the
open convex set $U\subseteq\mathbb{R}$. The following three statements are equivalent: (a) there exists a real-valued function $\theta$ such that $F(x) = \nabla\theta(x) \; \forall x \in U$; (b) the Jacobian matrix of $F(x)$ is symmetric $\forall x \in U$; (c) $F$ is integrable on $U$.
\end{theorem}
In terms of equilibrium problem \eqref{equilibrium_model}, $K=\set{O}\times\set{K}\times\Lambda_{+}$, vector $\textcolor{black}{x = [
p,
\boldsymbol{r},
d,
\boldsymbol{l},
\boldsymbol{\tilde{\lambda}}]^{\top}}$, and
$$
F^{\top} = [\begin{smallmatrix}
\nabla_{p}J^{\text{p}}(p,\boldsymbol{r}) &
\nabla_{\boldsymbol{r}}J^{\text{p}}(p,\boldsymbol{r}) &
\nabla_{d}J^{\text{c}}(d,\boldsymbol{l}) &
\nabla_{\boldsymbol{l}}J^{\text{c}}(d,\boldsymbol{l}) &
\nabla_{\boldsymbol{\tilde{\lambda}}}\boldsymbol{J}^{\text{ps}}(\boldsymbol{\tilde{\lambda}})
\end{smallmatrix}],
$$
where symbols in bold are properly dimensioned vectors.
\subsection{Proof}
\circled{\normalsize{1}} Existence. Recall that by definition $\set{O}$, $\set{K}$ and $\Lambda_{+}$ are convex and compact. The map $F$ is continuous as agents' objective functions are differentiable. Thus, the solution to equilibrium exists by Theorem \ref{theorem_ex}.

\circled{\normalsize{2}} Uniqueness. We rely on the symmetry principle that states that if Jacobian of $F$ is symmetric, there exists an equivalent optimization problem that solves $\mathrm{VI}$($K$,$F$).  \textcolor{black}{The Jacobian writes as:}
\textcolor{black}{
\begin{equation*}
\nabla_{x}F(x)=
\begingroup % keep the change local
\setlength\arraycolsep{0.5pt}
\begin{pmatrix*}[c]
\alpha p&       0&      0&      0&      -1 \\
0&  \alpha \boldsymbol{\pi}^{\text{p}\top} \boldsymbol{r}& 0&  0&      \boldsymbol{-1} \\ 
0& 0& \beta d & 0& 1 \\
0& 0& 0& \beta \boldsymbol{\pi}^{\text{c}\top} \boldsymbol{l} & \boldsymbol{1} \\
1&\boldsymbol{1}& -1& -\boldsymbol{1}& 0
\end{pmatrix*},
\endgroup
\end{equation*}
which includes a symmetric part with entries corresponding to the elements of variable set  $x^{\prime}=\{p,\boldsymbol{r},d,\boldsymbol{l}\}$. We further observe that $\nabla_{x^{\prime}}F(x^{\prime})$ is continuous in $x^{\prime}$, thus the conditions (b,c) of Theorem \ref{theorem_un} hold for the symmetric part, such that there exists a function $\theta(x^{\prime})$ given by
\begingroup
\allowdisplaybreaks
\begin{align*}
    &\theta(x^{\prime}) = \int_{0}^{1} F(x_{0}^{\prime} + t(x^{\prime}-x_{0}^{\prime}))^{\top}(x^{\prime}-x_{0}^{\prime})dt \\
    &\overset{\textcolor{black}{x_{0}^{\prime}\rightarrow0}}{=} \int_{0}^{1}
    \begin{pmatrix*}[l]
    t \alpha p \\
    t \alpha \boldsymbol{\pi}^{\text{p}\top} \boldsymbol{r}\\
    -\gamma + t \beta d\\
    - \gamma \boldsymbol{\pi}^{\text{c}} + t \beta \boldsymbol{\pi}^{\text{c}\top} \boldsymbol{l}
    \end{pmatrix*}^{\top}
    \begin{pmatrix}
    p \\
    \boldsymbol{r} \\
    d \\
    \boldsymbol{l}
    \end{pmatrix}dt \\
    &= [\alpha p^2 
    + \Sigma_{\omega}\pi_{\omega}^{\text{p}} \alpha r_{\omega}^2 
    + \beta d^2 
    + \Sigma_{\omega}\pi_{\omega}^{\text{c}} \beta l_{\omega}^2]
    \int_{0}^{1} t dt\\
    &- [\gamma d + \Sigma_{\omega}\pi_{\omega}^{\text{c}}l_{\omega}]
    \\
    &=
    \frac{1}{2}\alpha p^2 - [\gamma d - \frac{1}{2} \beta d^2] + \Sigma_{\omega}\pi_{\omega}^{\text{p}} \frac{1}{2}\alpha r_{\omega}^2  - \Sigma_{\omega}\pi_{\omega}^{\text{c}} [\gamma l_{\omega} - \frac{1}{2}\beta l_{\omega}^2].
\end{align*}
\endgroup
If we optimize $\theta(x^{\prime})$ subject to the stationarity conditions of the price-setting agent, we derive the following optimization:
}
\begingroup
\allowdisplaybreaks
\begin{subequations}\label{eq_opt}
\begin{align}
    \maximize {p,r_{\omega},d,l_{\omega}} & \big[u(d)-c(p)\big] + \sum_{\omega \in \Omega}\pi_{\omega}^{\text{c}}u(l_{\omega})-\sum_{\omega \in \Omega}\pi_{\omega}^{\text{p}}c(r_{\omega}),\\
    \st \quad & p + r_{\omega} + \xi_{\omega} - d - l_{\omega} \geq 0 : \tilde{\lambda}_{\omega}, \quad \forall \omega \in \Omega,\\
    & (p,r_{\omega}) \in \set{O}, \; (d,l_{\omega}) \in \set{K}, \quad \forall \omega \in \Omega,
\end{align}
\end{subequations}
\endgroup
whose stationarity conditions correspond to \textcolor{black}{those of equilibrium problem \eqref{equilibrium_model}}. We know that optimization \eqref{eq_opt} yields a unique solution due to strict concavity of objective function and convex and compact constraint set. Since the solution of \eqref{eq_opt} constitutes set $\mathrm{SOL}$($K$,$F$), the solution of original equilibrium problem \eqref{equilibrium_model} is also unique, as desired.

\vspace{-0.4cm}
\section{Proof of Proposition \ref{Prop_equiv}.}\label{Proof_Prop_1}
Since producer and consumer optimize over independent variables, we can optimize problems \eqref{profit_max} and \eqref{utility_max} jointly. If we constrain the joint problem by the optimality conditions of price-setter problem \eqref{price_setter}, we obtained the following optimization:
\begin{subequations}
\begingroup
\allowdisplaybreaks
\begin{align}
    \maximize{p,r_{\omega},d,l_{\omega}} & J^{\text{p}}(p,r_{\omega}) + J^{\text{c}}(d,l_{\omega}), \\
    \text{s.t.}\quad& (p,r_{\omega}) \in \set{O}, \; (d,l_{\omega}) \in \set{K}, \; \forall \omega \in \Omega,\\
    &0 \leq p + r_{\omega} + \xi_{\omega} - d - l_{\omega} \perp \tilde{\lambda}_{\omega} \geq 0, \; \forall \omega.
\end{align}
\endgroup
\end{subequations}
This is equivalent to the optimality condition of the centralized problem \eqref{stoch_problem} given that expectations over uncertain renewable production are the same.

\vspace{-0.4cm}
\section{Proof of Proposition \ref{Prop_price_setter_solution}}
The descent direction of the price-setter problem writes as
$$-\nabla_{\tilde{\lambda}_{\omega}}J_{\omega}^{\text{ps}}(\tilde{\lambda}_{\omega}) = p^{\nu} + r_{\omega}^{\nu} + \xi_{\omega} - d^{\nu} - l_{\omega}^{\nu}.$$
Then, the solution of the price-setter problem evolves along the decent direction with a suitable step size $\rho$ as follows:
$$\tilde{\lambda}_{\omega}^{\nu} = \tilde{\lambda}_{\omega}^{\nu-1} - \rho\nabla_{\tilde{\lambda}_{\omega}}J_{\omega}^{\text{ps}}(\tilde{\lambda}_{\omega}),$$
that is bounded from below by zero due to $\tilde{\lambda}_{\omega} \in \Lambda_{+}$.

\end{document}